
\newcount\secno
\newcount\prmno
\newif\ifnotfound
\newif\iffound

\def\namedef#1{\expandafter\def\csname #1\endcsname}
\def\nameuse#1{\csname #1\endcsname}

\long\def\ifundefined#1#2#3{\expandafter\ifx\csname
  #1\endcsname\relax#2\else#3\fi}
\def\hwrite#1#2{{\let\the=0\edef\next{\write#1{#2}}\next}}

\toksdef\ta=0 \toksdef\tb=2
\long\def\leftappenditem#1\to#2{\ta={\\{#1}}\tb=\expandafter{#2}%
                                \edef#2{\the\ta\the\tb}}
\long\def\rightappenditem#1\to#2{\ta={\\{#1}}\tb=\expandafter{#2}%
                                \edef#2{\the\tb\the\ta}}

\def\lop#1\to#2{\expandafter\lopoff#1\lopoff#1#2}
\long\def\lopoff\\#1#2\lopoff#3#4{\def#4{#1}\def#3{#2}}

\def\ismember#1\of#2{\foundfalse{\let\given=#1%
    \def\\##1{\def\next{##1}%
    \ifx\next\given{\global\foundtrue}\fi}#2}}

\def\section#1{\vskip1truecm
               \global\def\currenvir{section}
               \global\advance\secno by1\global\prmno=0
               {\bf \number\secno. {#1}}
               \smallskip}

\def\subsection{\smallskip
\global\def\currenvir{subsection}
                \global\advance\prmno by1
                \ind{(\number\secno.\number\prmno) }}
\def\subsec{\global\def\currenvir{subsection}
                \global\advance\prmno by1
                { (\number\secno.\number\prmno)\ }}

\def\proclaim#1{\global\advance\prmno by 1
                {\bf #1 \the\secno.\the\prmno$.-$ }}

\long\def\th#1 \enonce#2\endth{%
   \medbreak\proclaim{#1}{\it #2}\global\def\currenvir{th}\smallskip}

\def\bib#1{\rm #1}
\long\def\thr#1\bib#2\enonce#3\endth{%
\medbreak{\global\advance\prmno by 1\bf#1\the\secno.\the\prmno\ 
\bib{#2}$\!\cdot  -$ } {\it
#3}\global\def\currenvir{th}\smallskip}
\def\rem#1{\global\advance\prmno by 1
{\it #1} \the\secno.\the\prmno$.-$ }


\def\isinlabellist#1\of#2{\notfoundtrue%
   {\def\given{#1}%
    \def\\##1{\def\next{##1}%
    \lop\next\to\za\lop\next\to\zb%
    \ifx\za\given{\zb\global\notfoundfalse}\fi}#2}%
    \ifnotfound{\immediate\write16%
                 {Warning - [Page \the\pageno] {#1} No reference found}}%
                \fi}%
\def\ref#1{\ifx\labellist\empty{\immediate\write16
                 {Warning - No references found at all.}}
               \else{\isinlabellist{#1}\of\labellist}\fi}

\def\newlabel#1#2{\rightappenditem{\\{#1}\\{#2}}\to\labellist}
\def\labellist{}

\def\label#1{%
  \def\given{th}%
  \ifx\given\currenvir%
    {\hwrite\lbl{\string\newlabel{#1}{\number\secno.\number\prmno}}}\fi%
  \def\given{section}%
  \ifx\given\currenvir%
    {\hwrite\lbl{\string\newlabel{#1}{\number\secno}}}\fi%
  \def\given{subsection}%
  \ifx\given\currenvir%
    {\hwrite\lbl{\string\newlabel{#1}{\number\secno.\number\prmno}}}\fi%
  \def\given{subsubsection}%
  \ifx\given\currenvir%
  {\hwrite\lbl{\string%
    \newlabel{#1}{\number\secno.\number\subsecno.\number\subsubsecno}}}\fi
  \ignorespaces}

\newwrite\lbl

\def\openall{\openout\lbl=\jobname.lbl}

\newread\testfile
\def\lookatfile#1{\openin\testfile=\jobname.#1
    \ifeof\testfile{\immediate\openout\nameuse{#1}\jobname.#1
                    \write\nameuse{#1}{}
                    \immediate\closeout\nameuse{#1}}\fi%
    \immediate\closein\testfile}%

\def\begin{
\newlabel{corr}{1.1}
\newlabel{formules}{1.2}
\newlabel{dcorr}{1.3}
\newlabel{ch}{1.4}
\newlabel{ab}{1.5}
\newlabel{pol}{1.6}
\newlabel{polq}{2.6}
\newlabel{D}{3.1}
\newlabel{thm}{4.2}
\newlabel{lefdec}{5.3}
\newlabel{hardlef}{5.5}
\newlabel{h}{5.7}
\newlabel{wlef}{5.8}
\newlabel{conj}{5.9}
\newlabel{modconj}{5.10}
\newlabel{equiv}{5.11}
}


\magnification 1250
\pretolerance=500 \tolerance=1000  \brokenpenalty=5000
\mathcode`A="7041 \mathcode`B="7042 \mathcode`C="7043
\mathcode`D="7044 \mathcode`E="7045 \mathcode`F="7046
\mathcode`G="7047 \mathcode`H="7048 \mathcode`I="7049
\mathcode`J="704A \mathcode`K="704B \mathcode`L="704C
\mathcode`M="704D \mathcode`N="704E \mathcode`O="704F
\mathcode`P="7050 \mathcode`Q="7051 \mathcode`R="7052
\mathcode`S="7053 \mathcode`T="7054 \mathcode`U="7055
\mathcode`V="7056 \mathcode`W="7057 \mathcode`X="7058
\mathcode`Y="7059 \mathcode`Z="705A
\def\spacedmath#1{\def\packedmath##1${\bgroup\mathsurround =0pt##1\egroup$}
\mathsurround#1
\everymath={\packedmath}\everydisplay={\mathsurround=0pt}}
\def\nospacedmath{\mathsurround=0pt
\everymath={}\everydisplay={} } \spacedmath{2pt}
\def\qfl#1{\buildrel {#1}\over {\longrightarrow}}
\def\phfl#1#2{\normalbaselines{\baselineskip=0pt
\lineskip=10truept\lineskiplimit=1truept}\nospacedmath\smash {\mathop{\hbox to
8truemm{\rightarrowfill}}
\limits^{\scriptstyle#1}_{\scriptstyle#2}}}
\def\hfl#1#2{\normalbaselines{\baselineskip=0truept
\lineskip=10truept\lineskiplimit=1truept}\nospacedmath\smash{\mathop{\hbox to
12truemm{\rightarrowfill}}\limits^{\scriptstyle#1}_{\scriptstyle#2}}}
\def\diagram#1{\def\normalbaselines{\baselineskip=0truept
\lineskip=10truept\lineskiplimit=1truept}   \matrix{#1}}
\def\vfl#1#2{\llap{$\scriptstyle#1$}\left\downarrow\vbox to
6truemm{}\right.\rlap{$\scriptstyle#2$}}

\def\iso{\vbox{\hbox to .8cm{\hfill{$\scriptstyle\sim$}\hfill}
\nointerlineskip\hbox to .8cm{{\hfill$\longrightarrow $\hfill}} }}

\def\sdir_#1^#2{\mathrel{\mathop{\kern0pt\oplus}\limits_{#1}^{#2}}}
\def\pprod_#1^#2{\raise
2pt \hbox{$\mathrel{\scriptstyle\mathop{\kern0pt\prod}\limits_{#1}^{#2}}$}}

\font\eightrm=cmr8         \font\eighti=cmmi8
\font\eightsy=cmsy8        \font\eightbf=cmbx8
\font\eighttt=cmtt8        \font\eightit=cmti8
\font\eightsl=cmsl8        \font\sixrm=cmr6
\font\sixi=cmmi6           \font\sixsy=cmsy6
\font\sixbf=cmbx6\catcode`\@=11
\def\eightpoint{%
  \textfont0=\eightrm \scriptfont0=\sixrm \scriptscriptfont0=\fiverm
  \def\rm{\fam\z@\eightrm}%
  \textfont1=\eighti  \scriptfont1=\sixi  \scriptscriptfont1=\fivei
  \def\oldstyle{\fam\@ne\eighti}\let\old=\oldstyle
  \textfont2=\eightsy \scriptfont2=\sixsy \scriptscriptfont2=\fivesy
  \textfont\itfam=\eightit
  \def\it{\fam\itfam\eightit}%
  \textfont\slfam=\eightsl
  \def\sl{\fam\slfam\eightsl}%
  \textfont\bffam=\eightbf \scriptfont\bffam=\sixbf
  \scriptscriptfont\bffam=\fivebf
  \def\bf{\fam\bffam\eightbf}%
  \textfont\ttfam=\eighttt
  \def\tt{\fam\ttfam\eighttt}%
  \abovedisplayskip=9pt plus 3pt minus 9pt
  \belowdisplayskip=\abovedisplayskip
  \abovedisplayshortskip=0pt plus 3pt
  \belowdisplayshortskip=3pt plus 3pt 
  \smallskipamount=2pt plus 1pt minus 1pt
  \medskipamount=4pt plus 2pt minus 1pt
  \bigskipamount=9pt plus 3pt minus 3pt
  \normalbaselineskip=9pt
  \setbox\strutbox=\hbox{\vrule height7pt depth2pt width0pt}%
  \normalbaselines\rm}\catcode`\@=12

\newcount\noteno
\noteno=0
\def\up#1{\raise 1ex\hbox{\sevenrm#1}}
\def\note#1{\global\advance\noteno by1
\footnote{\parindent0.4cm\up{\number\noteno}\
}{\vtop{\eightpoint\baselineskip12pt\hsize15.5truecm\noindent
#1\medskip
}}\parindent 0cm}

\def\pc#1{\tenrm#1\sevenrm}
\def\tx{\kern-1.5pt -}
\def\cqfd{\kern 2truemm\unskip\penalty 500\vrule height 4pt depth 0pt width
4pt\medbreak} 
\def\virg{\raise
.4ex\hbox{,}}
\def\decale#1{\smallbreak\hskip 28pt\llap{#1}\kern 5pt}
\def\no{n\up{o}\kern 2pt}
\def\ind{\par\hskip 1truecm\relax}
\def\indp{\par\hskip 0.5truecm\relax}
\def\moins{\mathrel{\hbox{\vrule height 3pt depth -2pt width 6pt}}}
\def\rond{\kern 1,5pt\raise 0,5pt\hbox{${\scriptstyle\circ}$}\kern 2pt}
\def\End{\mathop{\rm End}\nolimits}

\def\Aut{\mathop{\rm Aut}\nolimits}

\frenchspacing
\input amssym.def
\input amssym
\vsize = 24truecm
\hsize = 16truecm
\hoffset = 0,4truecm
\voffset = -0,4truecm
\parindent=0cm
\baselineskip15pt
\overfullrule=0pt

\def\Q{{\bf Q}}
\def\Z{{\bf Z}}

\def\G{{\bf G}}
\def\slz{{SL_2(\Z)}}
\def\sp{{\bf SL}_2}
\def\slg{{\goth s}{\goth l}_2(\Q)}
\def\ch{\mathop{\rm ch}\nolimits}
\def\co{\mathop{\rm Corr}\nolimits}
\def\dco{\mathop{\rm Dcorr}\nolimits}
\def\bco{\mathop{\bf Corr}^*\nolimits}
\begin
\centerline{\bf  The action of $\sp$ on abelian varieties}\smallskip
 \centerline{Arnaud {\pc BEAUVILLE}} 
\vskip1.2cm

{\bf Introduction}
\smallskip
\ind The title is  somewhat paradoxical: we know that a linear group can only act trivially on an abelian variety. However we also know that there are not enough morphisms in algebraic geometry, a problem which may be fixed sometimes by considering {\it correspondences} between two varieties -- that is, algebraic cycles on their product, modulo rational equivalence. Our main result is the construction of a natural morphism of the algebraic group $\sp$ into the group $\co(A)^* $ of (invertible) self-correspondences of any polarized abelian variety $A$. As a consequence the group $\sp$ acts on the $\Q$\tx vector space $CH(A)$ parametrizing algebraic cycles (with rational coefficients) modulo rational equivalence, in such a way that this space decomposes as the direct sum of irreducible finite-dimensional representations.  This gives various results of Lefschetz type for the Chow group.
\ind This action of $\sp$ on $CH(A)$ is already known: it appears implicitely in the work of K\"unnemann [K], and explicitely in the unpublished thesis [P2]. But though it has been repeatedly used in  recent work on the subject ([P3], [P4], 
[Mo]), a detailed exposition does not seem to be available in the literature. The aim of this 
paper is to fill this gap, and also to explain the link with the action of $\slz$ on the derived category ${\bf D}(A)$ found by Mukai [M]. 

\ind In section 1 we recall some classical facts on  correspondences, mainly to fix our notations and conventions. In section 2 we explain how to deduce from Mukai's results a homomorphism  of $\slz$  into $\co(A)^* $, hence an action of $\slz$ onto $CH(A)$. 
 In sections 3 and 4 we show that these extend to  $\sp$, using a description 
of this algebraic group by generators and relations due to Demazure. In section 5 we 
deduce  some applications; the most interesting perhaps is a twisted version of the hard 
Lefschetz theorem for $CH(A)$: 
if $\theta \in CH^1(A)$ is an ample symmetric class, the multiplication map
$\times\,\, \theta ^{g-2p+s}: CH^p_{s}(A)\rightarrow CH^{g-p+s}_{s}(A)$ is an 
isomorphism  (the subscript $s$ refers to the decomposition of $CH(A)$ into eigensubspaces for the operators $n^*_A$, $n\in\Z$, see (\ref{ab}) below).

 {\medskip {\eightrm\baselineskip=12pt
\leftskip1cm\rightskip1cm\hskip0.8truecm  I am  grateful  to B. Fu for pointing out  Proposition \ref{equiv}, which led me to the Lefschetz type results discussed in section 5.\par}}\smallskip

\section{Correspondences}
\subsection We fix an algebraically closed field $k$. We denote by $A$ a smooth projective variety over $k$; from (\ref{ch}) on $A$ will be an abelian variety.


\ind As mentioned in the introduction, we will denote by $CH(A)$ the Chow group of algebraic cycles {\it with rational coefficients} on $A$ modulo rational equivalence\note{We could replace rational equivalence by any adequate equivalence relation, see [S].}. We briefly recall some basic facts about correspondences,  referring to [F] for a detailed treatment.
\ind A {\it correspondence} of $A$  is an element of $CH(A\times A)$. If $\alpha ,\beta $ are two  correspondences, we define their composition by $\beta \rond \alpha =(p_{13})_*(p_{12}^*\alpha  \cdot  p_{23}^*\beta )$, where $p_{ij}: A\times A\times A\rightarrow A\times A$ is the projection on the $i$\tx th and $j$\tx th factors. This defines an (associative) $\Q$\tx algebra structure on $CH(A\times A)$; this algebra is denoted $\co(A)$, and its subgroup of invertible elements by $\co(A)^*$. The unit element is the class $[\Delta _A]$ of the diagonal in $A\times A$. 
\ind To a class $\alpha $ in $\co(A)$ we associate a $\Q$\tx linear map 
$$\alpha _*:CH(A)\rightarrow CH(A)\quad\hbox{defined by }\ \alpha _*z=q_*(\alpha \cdot p^*z)\ ,$$ 
where $p$ and $q$ are the two projections from $A\times A$ to $A$. This gives a  $\Q$\tx algebra homomorphism
$\co(A)\rightarrow \End_{\Q}(CH(A))$, hence a group homomorphism $\co(A)^*\rightarrow \Aut_{\Q}(CH(A))$.

\label{corr}
\subsection We will need a few formulas satisfied by correspondences; they  can be found in (or follow easily from) [F],\S 16.
\ind{(\number\secno.\number\prmno.{\it a})} Let $\Delta :A\hookrightarrow A\times A$ be the diagonal morphism. We have
$$\Delta _*z \rond \alpha =\alpha \cdot q^*z\qquad \alpha \rond \Delta _*z  =\alpha \cdot p^*z\quad\hbox{for }\alpha \in\co(A)\,,\ z\in CH(A)\ ;$$taking $\alpha =\Delta _*x$ we see that {\it the map $\Delta _*:CH(A)\rightarrow \co(A)$ is a $\Q$\tx algebra homomorphism}.
\ind{(\number\secno.\number\prmno.{\it b})} Let $u$ be an endomorphism of $A$, $\Gamma _u$ the class of its graph, and $\Gamma '_u$ its transpose ($=(u,1_A)_*\Delta _A$). For $\alpha \in\co(A)$, we have
$$\Gamma _u\rond \alpha =(1,u)_*\alpha \ \ ,\ \  \Gamma' _u\rond \alpha =(1,u)^*\alpha \ \ ,\ \ 
\alpha \rond \Gamma _u=(u,1)^*\alpha \ \ ,\ \   \alpha\rond \Gamma' _u =(u,1)_*\alpha\ .$$
\label{formules}
\vskip-6pt
\subsection The above constructions work almost word for word replacing $CH(A)$ by the derived category ${\bf D}(A)$
of bounded  complexes of coherent sheaves on  $A$. A {\it derived correspondence} of $A$ is an object of ${\bf D}(A\times A)$. Using the notations of (\ref{corr}), we define the composition of two such objects $K,L$ as $L\rond K:= (p_{13})_*(p_{12}^*K\otimes p_{23}^*L)$. This defines an associative multiplication on the set of isomorphism classes of objects of  ${\bf D}(A\times A)$; we will denote this monoid by $\dco(A)$, and by $\dco(A)^*$ the subgroup  of invertible elements. Their unit is the 
sheaf ${\cal O}_{\Delta _A}$.
\ind As above we associate to $K\in\dco(A)$ the Fourier-Mukai transform $\Phi _K: $ ${\bf D}(A)\rightarrow {\bf D}(A)$, defined by $\Phi _K(\star)=q_*(p^*(\star)\otimes K)$. 
This defines   a 
group homomorphism
$$\Phi :\dco(A)^*\rightarrow \Aut({\bf D}(A))$$where $\Aut({\bf D}(A))$ is the group of isomorphism classes of auto-equivalences of the triangulated category  ${\bf D}(A)$.
By a celebrated theorem of Orlov [O], this map is {\it bijective}.\label{dcorr}

\subsection From now on we assume that $A$ is an abelian variety. In that case the  constructions of (\ref{corr}) and (\ref{dcorr}) are linked by the Chern character, which is a monoid homomorphism
$\ch :\dco(A)\rightarrow \co(A)$. We have a commutative diagram
$$\diagram{\dco(A)^* &\hfl{\Phi }{}&\Aut({\bf D}(A))\cr
\vfl{\ch}{}&&\vfl{}{\kappa }\cr
\co(A)^*&\hfl{*}{}&\Aut_{\Q}(CH(A))}$$where the map $\kappa $ is defined as follows: an automorphism of ${\bf D}(A)$ induces an automorphism of the K-theory group $K(A)$, hence a $\Q$\tx linear automorphism of $K(A)\otimes \Q$, hence 
a $\Q$\tx linear automorphism of $CH(A)$ via the isomorphism $\ch:K(A)\otimes \Q\iso CH(A)$.

\label{ch}

\subsection  For an abelian variety $A$, the unit element $[\Delta_A ]$ of  $\co(A)$ has a canonical decomposition    as a sum of orthogonal idempotents ([D-M], Theorem 3.1)
$$[\Delta_A ]=\sum_{i=0}^{2g} \pi _i$$
characterized by the property $(1,k)^*\pi _i=k^i\pi _i$. This decomposition induces a grading $CH(A)=\sdir_{s}^{}CH_s(A)$, with
 $$CH^p_s(A)= (\pi _{2p-s})_*(CH^p(A))=\{x\in CH^p(A)\ |\ n_A^*x=n^{2p-s}x \ \hbox{ for all }n\in \Z \}$$(see [B2]).

\label{ab}
\subsection Suppose now that $A$ has a polarization $\theta $, which we will view as a  {\it symmetric} (ample) element of $CH^1(A)$. The polarization defines an isogeny $\varphi : A\rightarrow \hat A$. The {\it  Poincar\'e bundle} ${\cal P}$ on $A\times A$  is by definition the pull back by $(1,\varphi )$ of the Poincar\'e bundle on $A\times \hat A$; it will play an important role in what follows.  

\ind Now there is a sign ambiguity in the definition of $\varphi $, hence of ${\cal P}$. Most authors, following 
Mumford in [Md], use the formula $\varphi (a)={\cal O}_A(\Theta -\Theta _a)$, where $
\Theta $ is a divisor  defining the polarization and $\Theta _a=\Theta +a$ denotes its 
translate by $a\in A$. This convention has some serious drawbacks. One of them is that if $A$ 
is  the Jacobian  of a curve $C$,  the natural 
map $\hat A={\bf Pic}^0(A)\rightarrow A={\bf Pic}^0(C)$ deduced from the embedding
 of $C$ in $A$ (defined up to translation) is the {\it opposite} of $\varphi ^{-1}$. More 
important for us,  it  leads to sign problems in the definition of  the action of (a covering 
of) $\slz$ on ${\bf D}(A)$. Because of these difficulties we will   
use the isomorphism $\varphi :A\rightarrow  \hat A$ defined by $\varphi (a)={\cal O}_A
(\Theta_a -\Theta)$. With this convention, by the see-saw theorem the class in $CH^1(A\times A)$ of the Poincar\'e bundle is $p^*\theta +q^*\theta -m^*\theta  $, where $m:A\times A\rightarrow A$ is the addition map $(a,b)\mapsto a+b$.
\label{pol}

\section{The homomorphism  $\slz\rightarrow \co(A)^*$}

\subsection In [M] Mukai observes that the derived category ${\bf D}(A)$ of a principally polarized abelian variety $(A,\theta )$ carries an 
action of $\slz$ ``up to shift". This is nicely elaborated in [P1] as an action of  a central extension of 
$\slz$ by $\Z$, the trefoil group  $\widetilde{SL}_2(\Z) $ (also known as the braid group on three strands). We will need only to describe this action in a
na\"{\i}ve sense, that is, as a group homomorphism 
of $\widetilde{SL}_2(\Z) $ in the group $\Aut({\bf D}(A))$ (\ref{dcorr}).

\ind  Recall that the group $\slz$  is generated by the elements
$$w=\pmatrix{0 &-1\cr 1&0}\qquad u=\pmatrix{1&1\cr 0&1}$$
with the relations  $w^2=(uw)^3$, $w^4=1$. The group    $\widetilde{SL}_2(\Z) $ is generated by two   
elements $\tilde u, \tilde w$ with the relation $\tilde w^2=(\tilde u\tilde w)^3$; the 
covering $\widetilde{SL}_2(\Z) \rightarrow \slz$ is obtained by mapping $\tilde u$ to $u$ and $\tilde w$ to 
$w$.
\subsection Let  ${\cal P}$ be the {\it Poincar\'e line bundle} on $A\times A$ (\ref{pol}). The functor $\Phi _{\cal P}$ is an autoequivalence of the category ${\bf D}(A)$ -- this is the original Fourier-Mukai functor [M]. We
choose a {\it symmetric} line bundle $L$ on $A$ with class $\theta $ and
 define an action of $\widetilde{SL}_2(\Z) $ on ${\bf D}(A)$ by mapping $\tilde u$ to the functor 
$\otimes L$ and
 $\tilde w$ to  $\Phi _{\cal P}$. Theorem 3.13 of [M] gives\note{Note 
that our  functor $\Phi _{\cal P}$ is Mukai's $R{\cal S}$  composed with $(-1_A)^*$.} $(\Phi _{\cal P} )^2=
(\otimes L\rond \Phi _{\cal P} )^3=(-1_A)^*[-g]$, so that we have indeed an action of 
$\widetilde{SL}_2(\Z)  $ on ${\bf D}(A)$, with the central element $z=
\tilde w^2=(\tilde u\tilde w)^3$ acting as $(-1_A)^*[-g]$; thus $z^2$ acts as the shift  $[-2g]$.

\subsection Observe that the functor  $\otimes L$ can be written $\Phi _{\Delta _*L}$. By (\ref{dcorr}) it follows that the homomorphism $\widetilde{SL}_2(\Z)  \rightarrow \Aut({\bf D}(A))$ factors through a group homomorphism
$$\widetilde{SL}_2(\Z)  \rightarrow \dco(A)^*$$
mapping $\tilde u$ to $\Delta _*L$ and $\tilde w$ to ${\cal P}$. 

\subsection We now consider the composition $\slz\rightarrow \dco(A)^*\qfl{\ch} \co(A)^*$, where $\ch$ is the Chern character (\ref{ch}).
Since $\ch E[p]=(-1)^p\ch E$,  this homomorphism maps $z^2$ to the unit element, hence factors as a homomorphism
$$  \slz \rightarrow \co(A)^*$$
which maps $u$ to  $\Delta _*e^\theta$, and $w$ to $e^\wp$, where $\wp$ is the class of ${\cal P}$ in $CH^1(A)$.

\subsection  The  argument extends with little change to the case of an arbitrary 
polarization. Let $A$ be a polarized abelian variety, of dimension $g$; we  denote by $\theta $ the 
unique symmetric element in $CH^1(A)$ representing the polarization,   and by $\wp=p^*\theta +q^*\theta -m^*\theta$  the  class  in $CH^1(A\times A)$ of the Poincar\'e bundle (\ref{pol}). The {\it degree} of $\theta $ is $d={\theta ^g\over g!}$.

  \th Proposition
 \enonce There is a {\rm (}unique{\rm )} group homomorphism $\slz\rightarrow \co(A)^*$ mapping $u$ to $\Delta _*e^\theta$ and $w$ to $d^{-1}e^{\wp}$.
 
 \endth\label{polq}
{\it Proof} :  We choose an isogeny $\pi$ of $A$ onto an abelian variety $A_0$ with a 
principal polarization $\theta _0$ such that $\theta =\pi ^*\theta _0$ ([Md], \S 23, Cor. 1 
of Thm. 4). One checks readily that the $\Q$\tx linear isomorphism   $d^{-1}(\pi ,\pi )^*:CH(A_0\times A_0)\rightarrow CH(A\times A)$ is compatible with the composition of correspondences, thus induces an isomorphism of algebras
$\co(A_0)\iso \co(A)$. 
Let $\Delta _0$ denote the diagonal morphism of $A_0$. We have $(\pi ,\pi )\rond \Delta=\Delta _{0}\rond \pi $ and $f\rond (\pi ,\pi )=\pi \rond f$ if $f=p,q$ or $m$. From this one easily gets

$$(\pi ,\pi )^*\Delta _*e^{\theta _0}=d\Delta _*e^\theta \ \hbox{ and }\ (\pi ,\pi )^*e^{p^*\theta _0+q^*\theta _0-m^*\theta _0}=e^{p^*\theta +q^*\theta -m^*\theta }\ ,$$hence the result.\cqfd

\section{Extension to  $\sp$}
\ind We will now show that the homomorphism  $\slz \rightarrow \co(A)^*$ extends to  the 
algebraic group $\sp$  over $\Q$.  
The essential tool is the description of $\sp$ by generators and 
relations given (in a much more general set-up) in [D], Theorem 6.2.  
\ind We denote by $B$  the upper triangular Borel subgroup of $\sp$. We still denote by $w$ 
and $u$ the elements $\pmatrix{0&-1\cr 1&0}$ and $\pmatrix{1&0\cr 1&1}$ of $SL_2
(\Q)$. By a $\Q$\tx group we mean a sheaf of groups over ${\rm Spec}(\Q)$ for the fppf topology.

\th Proposition
\enonce Let  $H$ be a $\Q$\tx group. Suppose given a morphism of $\Q$\tx groups $\beta :B
\rightarrow H$ and an element $h\in H(\Q)$. Assume that:

\indp {\rm (i)} $h \beta(t) h^{-1}=\beta(t^{-1})$ for $t$ in the maximal torus of $B$;
\indp{\rm  (ii)} $h^2=(\beta(u) h)^3=\beta (-I)$ in $H(\Q)$.

Then there is a (unique) morphism of $\Q$\tx groups $\varphi :\sp\rightarrow H$ 
extending $\beta$ and mapping $w$ to $h$.

\endth\label{D}

{\it Proof} : This is   [D], Theorem 6.2, in the case $S={\rm Spec}(\Q)$, $G=\sp$ (note that our element $w$ is the opposite of
 the one in {\it loc. cit.}).\cqfd
\ind Observe that (ii) can be rephrased by saying that there is a homomorphism $\slz
\rightarrow H(\Q)$ which maps $w$ to $h$ and coincides with $\beta $ on $B(\Z)$.

\subsection If $C$ is a $\Q$\tx algebra, the functor $R\mapsto (C\otimes _\Q R)^*$ is a $\Q$\tx group; its Lie algebra is $C$, endowed with the bracket $[x,y]=xy-yx$. We will denote by $\bco(A)$ the $\Q$\tx group obtained from the $\Q$\tx algebra $\co(A)$ in this way.

\ind Let $\delta :A\times A\rightarrow A$ be the difference map $(a,b)\mapsto b-a$. For $n\in\Z$ we denote by $\Gamma _n\in \co(A)$ the graph of the multiplication by $n$, and by $\Gamma '_n$ its transpose.

\th Theorem
\enonce Let $A$ be an abelian variety, of dimension $g$, with a polarization $\theta $ of degree $d$. There 
is a morphism  of $\Q$\tx groups $\varphi :\sp \rightarrow \bco(A)$ such that, for $n\in\Z\moins\{0\}$, $a\in\Q$:
$$\nospacedmath\displaylines{\varphi \Bigl(\pmatrix{n &0\cr 0&n^{-1}} \Bigr)=n^{-g}  \Gamma' _n
\quad,\quad 
\varphi \Bigl(\pmatrix{0 &-1\cr 1&0}\Bigr)=d^{-1}e^{\wp} \ , \cr
\varphi \Bigl(\pmatrix{1 &a\cr 0&1}\Bigr)=\Delta _*e^{a\theta}   \quad,\quad
\varphi \Bigl(\pmatrix{1 &0\cr a&1}\Bigr)=d^{-1}a^g\,e^{\delta ^*\theta/a}\ .
}$$
\ind The corresponding Lie algebra homomorphism $L\varphi :\slg\rightarrow \co(A)$ is given, in the standard basis $(X,Y,H)$ of $\slg$, by:
$$L\varphi (X)=\Delta _*\theta  \quad,\quad  L\varphi (Y)=  {\delta ^*\theta ^{g-1}\over d(g-1)!}  \quad,\quad  L\varphi (H)=\sum_i(i-g)\pi _i\ .$$
\endth
{\it Proof} : We  apply Proposition \ref{D} with $H=\bco(A)$. 
To define $\beta $ we write $B$ as a semi-direct product $\G_a\rtimes \G_m$. We define $\alpha:\G_a\rightarrow \bco(A)$ by $a\mapsto \Delta _*e^{a\theta }$; this is a morphism of $\Q$\tx groups by (\ref{formules}.{\it a}).

\ind We define $\tau:\G_m\rightarrow \bco(A)$ by
$\tau(t)=t^{-g}\sum_i t^{i}\pi _i$, where $(\pi _i)$ is the family of orthogonal idempotents considered in (\ref{ab}).
This is a morphism of $\Q$\tx groups; for $t\in\Z$ we have $\tau(t)=t^{-g} (1,t)^*\Delta_A =t^{-g}\Gamma' _t$ and $\tau (t^{-1})=t^g(\Gamma '_t)^{-1}=t^{-g}\Gamma _t$. 

\ind To ensure that $\beta=(\alpha,\tau):\G_a\rtimes \G_m\rightarrow \bco(A)$ is a morphism of groups, we must check the commutation relation  $\tau(t)\,\alpha(a)\,\tau(t)^{-1}=\alpha (t^2a)$. This relation is polynomial in $t$, so it suffices to check it for $t=s^{-1}$ with $s\in\Z$. In that case we have by (\ref{formules}.{\it b}):
$$\tau(t)\,\alpha(a)\,\tau(t)^{-1}=s^{-2g} \Gamma _s\rond \Delta _*e^{a\theta }\rond \Gamma' _s=s^{-2g}(s,s)_*\Delta_*e^{a\theta }=s^{-2g}\Delta _*s_*e^{a\theta }\ ;$$
since $s_*\theta ^p=s^{2g-2p}\,\theta ^p$, this gives $\tau(t)\,\alpha(a)\,\tau(t)^{-1}=\Delta _*e^{t^2a\theta }=\alpha (t^2a)$ as required.

\ind We take for $h$ the element $d^{-1}e^{\wp}$ of $\co(A)^*$. Condition (ii) is then satisfied because $h$ and $\beta(u)$ are the images of $w$ and $u$ by the homomorphism $\slz\rightarrow \co(A)^*$ (Prop. \ref{polq}).
\ind Condition (i) can be written $\tau (t)h\tau(t)=h$. Again it suffices  to check the equality
$\Gamma' _t\rond e^{\wp}\rond \Gamma' _t=t^{2g}e^{\wp}$ for $t\in\Z$.
The Poincar\'e bundle ${\cal P}$ satisfies $(t,1)^*{\cal P}=$ $(1,t)^*{\cal P}={\cal P}^t$ for every $t\in\Z$. Therefore 
$(t,1)^*\wp=(1,t)^*\wp$, and 
$$\Gamma' _t\rond e^{\wp}\rond \Gamma' _t = (t,1)_*(t,1)^*e^{\wp}=t^{2g}e^{\wp}\ . $$
\ind This proves the existence of $\varphi $ satisfying the three first formulas stated. To prove the fourth one, put
 $v=\pmatrix{1 & 0\cr 1 & 1}$; we have $v= uwu$. Thus, using (\ref{formules}.{\it a}),
$$\varphi (v)=\Delta _*e^{\theta} \rond e^\wp\rond \Delta _* e^{\theta} =  e^{\wp+p^*\theta +q^*\theta  }\ .$$
Since $\wp=p^*\theta +q^*\theta -m^*\theta$ and $m^*\theta +\delta  ^* \theta =2p^*\theta +2q^*\theta$ by the seesaw theorem, this gives $\varphi (v)=e^{\delta ^*\theta }$.


\ind Now we use the equality $\pmatrix{t &0\cr 0&t^{-1}}\pmatrix{1 &0\cr 1&1}\pmatrix{t^{-1} &0\cr 0&t}=\pmatrix{1 &0\cr t^{-2}&1} $ to get, for $t\in\Z$,

$$\varphi \Bigl(\pmatrix{1 &0\cr t^{-2}&1}\Bigr)= t ^{-2g}\ \Gamma' _t\rond e^{\delta ^*\theta } \rond \Gamma _t$$
We have by (\ref{formules}.{\it b})
$$\Gamma' _t\rond e^{\delta ^*\theta} \rond \Gamma _t=(t,t)^*e^{\delta ^*\theta} =e^{\delta ^*t^*\theta }=e^{t^2\delta ^*\theta }\ ;$$putting $a=t^{-2}$ gives  $\varphi \Bigl(\pmatrix{1 &0\cr a&1}\Bigr)=d^{-1}a^g\,e^{\delta ^*\theta/a}$ for $a=t^{-2}$ with $t\in\Z$, hence as usual for all $a$.
\ind The value of $L\varphi $ follows from these formulas by differentiation.\cqfd

\section{The action of $\sp$ on ${\bf CH(A)}$}

\subsection 
Let $V$ be a $\Q$\tx vector space; the functor $R\mapsto  GL(V\otimes R)$ on the category of (commutative) $\Q$\tx algebras is a $\Q$\tx group, that we will denote by ${\bf GL}(V)$. If $G$ is an algebraic group  over $\Q$, we define a representation of $G$ on $V$ as a morphism of $\Q$\tx groups $G\rightarrow {\bf GL}(V)$.

\ind Recall that 
the {\it Pontryagin product}  of two elements $\alpha ,\beta $ of $CH(A)$  is defined by $\alpha *\beta :=m_*(p^*\alpha \cdot q^*\beta )$.
\ind We will denote by ${\cal F}$ the  $\Q$\tx linear automorphism $d^{-1}(e^\wp)_*$ of $CH(A)$; this is the {\it Fourier transform} for Chow groups, see [B1] and [B2]. 

\th Theorem 
\enonce Let $A$ be an abelian variety, with a polarization $\theta $ of degree $d$. There 
is a representation of $\sp$ on $CH(A)$, which is a direct sum of finite-dimensional 
representations, such that, for $n\in\Z\moins\{0\}$, $a\in\Q$, $z\in CH(A):$
$$\nospacedmath\displaylines{\pmatrix{n &0\cr 0&n^{-1}}\!\cdot  z=n^{-g}n^*z  
\quad,\quad   
\pmatrix{0 &-1\cr 1&0}\!\cdot  z={\cal F} (z) \ , \cr
\pmatrix{1 &a\cr 0&1}\!\cdot  z=e^{a\theta}z   \quad,\quad
\pmatrix{1 &0\cr a&1}\!\cdot  z=d^{-1}a^ge^{\theta/a}*z\ .
}$$
\ind The corresponding action of the Lie algebra $\slg$ is given by:
$$Xz=\theta z\quad,\quad Yz=d^{-1}{\theta ^{g-1}\over (g-1)!}\, *\, z \quad,\quad  Hz=
(2p-g-s) z\quad {\rm for}\ \ z\in CH^p_s(A)\ .$$
\endth\label{thm}

{\it Proof} :  The homomorphism $\co(A)\rightarrow \End_{\Q}(CH(A))$ (\ref{corr}) defines a morphism of $\Q$\tx groups
$\bco(A)\rightarrow {\bf GL}(CH(A))$, hence by composition  with $\varphi  $  a representation of $\sp$ on $CH(A)$.  By definition, $g\cdot z=\varphi (g)_*z$ for $g\in\sp$, $z\in CH(A)$. Thus
\medskip

$\pmatrix{n &0\cr 0&n^{-1}}\!\cdot  z=n^{-g}(\Gamma _n)_*z=n^{-g}n^*z$,
$\quad\pmatrix{0 &-1\cr 1&0}\!\cdot  z=d^{-1}(e^\wp)_*z \buildrel {\scriptscriptstyle\rm def}\over {=}{\cal F} (z)$,
\medskip

$\pmatrix{1 &a\cr 0&1}\!\cdot  z=(\Delta _*e^{a\theta}z)_*= e^{a\theta}z$,
$\quad\pmatrix{1 &0\cr a&1}\!\cdot  z=d^{-1}a^g\,(e^{\delta ^*\theta/a})_*z$.
\medskip

Let $\sigma $ be the automorphism of $A\times A$ defined by $\sigma (a,b)=(b,a+b)$. We have 

$q\rond \sigma =m$, $p\rond \sigma =q$, $ \delta \rond \sigma =p$, hence  
$$(e^{\delta ^*\theta/a })_*z=q_*\sigma _*\sigma ^*(e^{\delta ^*\theta/a }\cdot p^*z)=m_*(p^*e^{\theta/a }\cdot q^*z)=e^{\theta/a} * z \ .$$

 In particular, $H$ is diagonalizable and $X,Y$ are nilpotent; this is enough to 
imply that $V$ is a direct sum of finite-dimensional representations  ([Bo], ch. 8, \S 1, exerc. 4).\cqfd

\section{Application to Lefschetz type results}
\subsection  In this section we will apply the well-known structure of finite-dimensional representations of $\sp$. We say that an element $z\in CH_s^p(A)$ is {\it primitive} if $\theta ^{g-1}
*z=0$. The primitive elements are exactly the lowest weight elements for the action of $
\sp$ on $CH(A)$. 
\ind Let $z\in CH_s^p(A)$ be a primitive element. The subspace of $CH(A)$ spanned (over $\Q$) by the $\theta ^qz$  is an irreducible representation of $\sp$; it is identified with the space of polynomials in one variable of degree $\leq g+s-2p$ (with the standard action) by the map $P\mapsto P(\theta )z$. This gives an explicit description of the action of $\sp$; in particular:

\th Proposition
\enonce If $z\in CH_s^p(A)$ is primitive, we have $g+s-2p\geq 0$, and $(z,\theta z,\ldots 
, \theta ^{g+s-2p}z)$ is a basis of an irreducible  subrepresentation of $CH(A)$. The vector space $CH(A)$ is a direct sum of subrepresentations of this type.\cqfd
\endth
\th Corollary
\enonce Let $P_s^p\i CH^p_s(A)$ be the subspace of primitive elements. Then $CH^p_s
(A)=\sdir_{q\leq p}^{}\theta ^{p-q}\,P_s^q$.\cqfd
\endth\label{lefdec}
\ind Since the Fourier automorphism ${\cal F}$  of $CH(A)$ is given by the action of $w$, we have:

\th Corollary
\enonce Let $z\in CH^p_s(A)$ be a primitive element, and let $q\leq g+s-2p$. Then $
\displaystyle {\cal F} \bigl({\theta ^q\over q!} z\bigr)= {(-\theta) ^{r}\over r!}z$, with $r=g+s-2p-q$.\cqfd
\endth

\th Proposition
\enonce  The multiplication map $\times\, \theta ^{ q-p} : CH^p_s(A) \rightarrow 
CH^q_s(A)$ $(q \geq p)$
is injective for $p + q  \leq  g + s$ and surjective for $p + q \geq   g + s $. In particular, it is bijective for $p+q=g+s$.
\endth\label{hardlef}
{\it Proof} : Assume  $p + q  \leq  g + s$; let $z\in CH^p_{s}(A)$ with $\theta ^{q-p}z=0$. Using Cor. \ref
{lefdec} we write 
$z=\sum_{r\leq p}\theta ^{p-r}z_r$, with $z_r\in P_s^r$; we have $\theta ^{q-r}z_r=0$ 
for each $r$. Since $q-r\leq g+s-2r$  this implies $z_r=0$ for each $r$, hence $z=0$.
\ind Assume $p + q  \geq  g + s$. To prove the surjectivity of $\times\, \theta ^{ q-p}$ it suffices, by  Cor. \ref{lefdec}, 
to prove that each nonzero element $\theta ^{q-r}z_r$, with $z_r\in P_s^r$, lies in the 
image. But since 
$\theta ^{q-r}z_r\neq 0$ we have $q-r\leq g+s-2r$, hence $q+r\leq g+s\leq p+q$ and finally $r\leq p$. Therefore
 $\theta ^{q-r}z_r=\theta ^{q-p}(\theta ^{p-r}z_r)$.\cqfd

\subsection In what follows we consider the filtration of $CH(A)$ associated to the grading 
(\ref{ab}), that is, $F^sCH^p(A) :=  \sdir_{t\geq s}^{}
CH^p_t(A)$. Conjecturally this is the Bloch-Beilinson filtration of $CH(A)$, see [Mu].
\th Corollary 
\enonce 
 Let $h \in CH^1(A)$ be an ample class. The multiplication map
$$\times\, h^{q-p} : F^{p+q-g}CH^p(A) \rightarrow  F^{p+q-g}CH^q(A)$$ is injective.
\endth\label{h}
{\it Proof} : Let $\theta $  be the component of $h$ in $CH^1_0
(A)$; it is ample and symmetric. The map induced by 
$\times\, h^{q-p}$ on the associated graded spaces is
$$\times\, \theta ^{q-p}: \sdir_{s\geq p+q-g}^{}CH^p_s(A)\longrightarrow  \sdir_{s\geq p+q-g}^{}
CH^q_s(A)\ ,$$
which is injective by Proposition \ref{hardlef}; therefore $\times\, h^{q-p}$ is injective.\cqfd
\th Corollary
\enonce Let $D$ be an ample divisor in $A $. The restriction map \break
$F^{2p+1-g}CH^p(A) \rightarrow  CH^p(D)$ is injective.
\endth\label{wlef}
\ind Here $CH^p(D)$ is the Chow group as defined in [F], chap. 2: $CH^p(D)=A_{g-1-p}
(D)$ in the notation of [F].

{\it Proof} : Let $i$ be the natural injection of $D$ in $A $, and let $h$ be the class of $D
$
in $CH^1(A) $. Let $z \in F^{2p+1-g}CH^p(A)$ such that $i^* z = 0 $. Then $h\cdot z = i_* 
i^* z = 0 $, so
$z = 0$ by Cor. \ref{h}.\cqfd

\subsection\label{conj} At this point we  recall the {\it vanishing conjecture} of [B2]:
$$CH_s(A)=0\quad{\rm for}\ s<0\ .\eqno{(\ref{conj})}$$
\ind This implies $F^0CH(A)=CH(A)$, hence:
\th Proposition
\enonce Assume that $CH_s(A)=0$ for $s<0$. Then:
\ind {\rm a)} If $h\in CH^1(A)$ is an ample class, the  multiplication map \break
$\times\,\, h^{q-p} : $ $CH^p(A) \rightarrow  CH^{q}(A)$ is injective provided $p+q\leq g$.
\ind {\rm b)} If $D$ is an ample divisor in $A$, the restriction map
$CH^p(A) \rightarrow  CH^p(D)$ is injective for $p\leq {1\over 2}(g-1)$.\cqfd
\endth\label{modconj}

\ind We have actually a more precise result,
which has been shown to me by B. Fu (with a different proof, see [Fu]):
\th Proposition 
\enonce Let  $h \in CH^1(A)$ be an ample class. The multiplication map
$\times\, h^{g-2p} : CH^p(A) \rightarrow  CH^{g-p}(A)$ is injective for all $p\leq {g\over 2}$ if and only if conjecture {\rm (\ref{conj})} holds.
\endth\label{equiv}
{\it Proof} : It remains to prove that if $CH_s(A)$ is nonzero for some $s<0$, all multiplication maps cannot be injective.  Let $\theta $  be the component of $h$ in $CH^1_0(A)$. By Cor. \ref{lefdec} there exists a nonzero primitive class $z\in CH^p_s(A)$ for some integer $p\leq {1\over 2}(g+s)$; we have $g-2p>g+s-2p$ and therefore $\theta ^{g-2p}z= 0$. The class $h$ is equal to $T_a^*\theta $ for some element $a\in A$, where $T_a$ denotes the translation $x\mapsto x+a$. We have $T_a^*z\equiv z$ mod.$F^{s+1}CH^p(A)$, hence $T_a^*z\neq 0$, and $h^{g-2p}(T_a^*z)=0$, which proves our assertion.\cqfd

\ind An  interesting feature  of Proposition \ref{modconj} is that it makes sense for any smooth projective variety $A$; the same is true of Cor. \ref{h}, taking for $(F^sCH(A))_{s\geq 0}$ the (conjectural) Bloch-Beilinson filtration of $CH(A)$. These conjectures are thoroughly discussed in [Fu], where it is shown that they would follow from a weak version of the Beilinson conjectures.

\vskip2cm
\font\cc=cmcsc10

\centerline{ REFERENCES} \vglue15pt\baselineskip12.8pt
\def\num#1#2#3{\smallskip\item{\hbox to\parindent{\enskip [#1]\hfill}}{\cc #2: }{\sl #3}}
\parindent=1.3cm 

\num{B1}{A. Beauville}{Quelques remarques sur la transformation de Fourier dans l'anneau de Chow d'une vari\'et\'e ab\'elienne}.  Algebraic geometry (Tokyo/\allowbreak Kyoto, 1982), 238--260, LNM {\bf 1016}, Springer, Berlin, 1983.
\num{B2}{A. Beauville}{Sur l'anneau de Chow d'une vari\'et\'e ab\'elienne}. Math. Ann. {\bf 273} (1986), no. 4, 647--651.
\num{Bo}{N. Bourbaki}{Groupes et alg\`ebres de Lie}, Chap.\ VIII. 
Hermann, Paris, 1975.

\num{D}{M. Demazure}{Groupes r\'eductifs de rang semi-simple}$\,1$. SGA 3, vol. III, Expos\'e 20. LNM {\bf 153}, Springer, Berlin, 1970.
\num{D-M}{C. Deninger, J. Murre}
{Motivic decomposition of abelian schemes and the Fourier transform}. 
J. Reine Angew. Math. {\bf 422} (1991), 201--219. 

\num{F}{W. Fulton}{Intersection theory}. Ergebnisse der Math. {\bf  2}. Springer-Verlag, Berlin, 1984.
\num{Fu}{B. Fu}{Remarks on
Hard Lefschetz conjectures on Chow groups}. Preprint.
\num{K}{K. K\"unnemann} {A Lefschetz decomposition for Chow motives of abelian schemes}. Invent. Math. {\bf 113} (1993), no. 1, 85--102.
\num{Mo}{B. Moonen}{Relations between tautological cycles on Jacobians}. Preprint arXiv:0706.3478. To appear in Comm. Math. Helvetici. 
\num{M}{S. Mukai}{Duality between $D(X)$ and $D(\hat X)$ with its application to Picard sheaves}. Nagoya Math. J. {\bf 81} (1981), 153--175.

\num{Md}{D. Mumford}{Abelian varieties}. Oxford University Press, London (1970).
\num{Mu}{J. P. Murre}
{On a conjectural filtration on the Chow groups of an algebraic variety, I. The general conjectures and some examples}. Indag. Math. (N.S.) {\bf 4} (1993), no. 2, 177--188. 
\num{O}{D. Orlov}{Equivalences of derived categories and $K3$ surfaces}. Algebraic geometry, 7. J. Math. Sci. (New York) {\bf 84} (1997), no. 5, 1361--1381.
\num{P1}{A. Polishchuk}{A remark on the Fourier-Mukai transform}. Math. Res. Lett. {\bf 2} (1995), no. 2, 193--202. 
\num{P2}{A. Polishchuk}{Biextension,Weil representation on derived categories, and theta functions}.
Ph. D. Thesis, Harvard Univ., 1996.
\num{P3}{A. Polishchuk}{Universal algebraic equivalences between tautological cycles on Jacobians of curves}. Math. Z. {\bf 251} (2005), no. 4, 875--897. 
\num{P4}{A. Polishchuk}{Lie symmetries of the Chow group of a Jacobian and the tautological subring}. J. Algebraic Geom. {\bf 16} (2007), no. 3, 459--476.
\num{S}{P. Samuel}{Relations d'\'equivalence en g\'eom\'etrie alg\'ebrique}. Proc. ICM 1958,  pp. 470--487; CUP, New York (1960).
\vskip1cm
\def\pc#1{\eightrm#1\sixrm}
\hfill\vtop{\eightrm\hbox to 5cm{\hfill Arnaud {\pc BEAUVILLE}\hfill}
 \hbox to 5cm{\hfill Institut Universitaire de France\hfill}\vskip-2pt
\hbox to 5cm{\hfill \&\hfill}\vskip-2pt
 \hbox to 5cm{\hfill Laboratoire J.-A. Dieudonn\'e\hfill}
 \hbox to 5cm{\sixrm\hfill UMR 6621 du CNRS\hfill}
\hbox to 5cm{\hfill {\pc UNIVERSIT\'E DE}  {\pc NICE}\hfill}
\hbox to 5cm{\hfill  Parc Valrose\hfill}
\hbox to 5cm{\hfill F-06108 {\pc NICE} cedex 2\hfill}}
\end